\documentclass[reqno,12 pt]{amsart}
\usepackage{latexsym,amssymb,amsmath,amsthm,amscd,amsthm,amsfonts,breqn,amssymb,enumerate,extarrows}
\usepackage{a4wide}
\usepackage{bbm}
\usepackage[all, knot]{xy}
\xyoption{arc}
\usepackage{mathrsfs}
\usepackage[ linkcolor=blue,
            anchorcolor=green,
            urlcolor=orange,
            citecolor=green
                         ]{hyperref}

\DeclareUnicodeCharacter{FFFD}{!!!FIX ME!!!} 
\DeclareFontFamily{U}{mathx}{}
\DeclareFontShape{U}{mathx}{m}{n}{ <-> mathx10 }{}
\DeclareSymbolFont{mathx}{U}{mathx}{m}{n}
\DeclareFontSubstitution{U}{mathx}{m}{n}

\DeclareMathAccent{\widecheck}{0}{mathx}{"71}

\DeclareUnicodeCharacter{FFFD}{!!!FIX ME!!!} 
\DeclareUnicodeCharacter{FF0C}{!!!FIX ME!!!} 
\usepackage{amsmath}
\allowdisplaybreaks[4]

\usepackage{amsmath}

\allowdisplaybreaks[4]

\makeindex

\DeclareFontFamily{U}{mathx}{}
\DeclareFontShape{U}{mathx}{m}{n}{ <-> mathx10 }{}
\DeclareSymbolFont{mathx}{U}{mathx}{m}{n}
\DeclareFontSubstitution{U}{mathx}{m}{n}

\DeclareMathAccent{\widecheck}{0}{mathx}{"71}

\newcommand{\ud}{\mathrm{d}}

\newcommand{\bee}{\begin{eqnarray}}
\newcommand{\ene}{\end{eqnarray}}
\newcommand{\bea}{\begin{eqnarray*}}
\newcommand{\ena}{\end{eqnarray*}}

\DeclareUnicodeCharacter{FFFD}{!!!FIX ME!!!} 
\DeclareUnicodeCharacter{0308}{!!!FIX ME!!!}  

\DeclareUnicodeCharacter{FB01}{!!!FIX ME!!!}  

\newcommand*{\C}{\ensuremath{\mathbb{C}}}
\newcommand*{\R}{\ensuremath{\mathcal{R}}}

\newcommand*{\Z}{\ensuremath{\mathbb{Z}}}
\newcommand*{\T}{\ensuremath{\mathcal{T}}}

\newcommand*{\GL}{\ensuremath{{\rm{GL}}}}
\newcommand*{\A}{\ensuremath{A_\Pi(n,1,1)}}
\newcommand*{\Kl}{\ensuremath{{\rm{Kl}_3}}}
\newcommand*{\TT}{\ensuremath{\mathfrak{t}}}
\newcommand*{\ST}{\ensuremath{\mathfrak{s}}}

\newcommand{\V}[1]{\ensuremath{V\lf (#1\ri)  }}

\newcommand{\U}[1]{\ensuremath{U\lf (#1\ri)  }}
\newcommand{\e}[1]{\ensuremath{e{\left( #1\right)  }}}
\newcommand{\fone}[1]{\ensuremath{\frac{1}{#1}}}

\newcommand{\lf}{\bigg}
\newcommand{\ri}{\bigg}
\newcommand{\f}{\frac}

\newcommand{\smu}{\hskip 0.5em \sideset{_{}^{}}{^{\ast}_{}}\sum\limits}

\theoremstyle{plain}
\newtheorem{theo}{Theorem}[section]
\newtheorem{lem}[theo]{Lemma}
\newtheorem{coro}[theo]{Corollary}

\theoremstyle{definition}

\makeatletter
\@addtoreset{equation}{section}
\makeatother

\numberwithin{equation}{section}

\DeclareUnicodeCharacter{FB03}{ffi}

\begin{document}

\title[]{\textbf {Hybrid bounds for ${\GL}(4)\times \text{\GL}(1)$ twisted $L$-functions}}



\author[F. Hou]{Fei Hou}
\address{School of Sciences, Xi'an University of Technology,
Xi'an 710054, China} \email{fhou@xaut.edu.cn}

\date{}
\maketitle

\begin{abstract}
\small {Let $P,M$ be two primes such that $(P,M)=1$. Let $\Pi$ be a normalized Hecke-Maa\ss\ form on $\text{GL(4)}$ of level $P$, and $\chi$ a primitive Dirichlet character modulo $M$. In this paper, we study the hybrid subconvexity problem for $L(s, \Pi\otimes \chi)$ simultaneously in the level and conductor aspects. Among other things, we prove a hybrid subconvex bound, so long as $M^{1/5}<P<M^{2/5}$. 
}
\end{abstract}
\section{Introduction}
Let $P,M$ be two primes such that $(P,M)=1$. Let $\chi $ be a primitive Dirichlet character modulo $M$. There is a very long history of understanding the central values of the $L$-functions. It might be traced back to the pioneering work of Burgess \cite{Bu}, who gave the well-known estimate
\[L(1/2,\chi)\ll    M^{{3}/{16}+\varepsilon} \]  in the year of 1962. Here, the exponent $3/16$ is nowadays known as the \emph{Burgess' quality}. In 1978, Heath-Brown \cite{HB}, for the first time, showed the hybrid subconvex bound for $L(s,\chi)$ simultaneously in the $s$ and conductor aspects. It is remarkable that, just recently Petrow and Young \cite{PY1} proved a Lindel\"{o}f-on-average upper bound for the fourth moment of Dirichlet $L$-functions along cosets, whereby to obtain the Weyl-strength hybrid bounds for all Dirichlet $L$-functions; this, in turn, improves on the result of Heath-Brown after roughly four decades.

Let $f$ be a Hecek newform on $\GL(2)$. In 2008, Blomer and Harcos \cite{BH} firstly studied hybrid subconvexity problem for ${\GL}(2)\times \GL(1) $ $L$-function $L(s,f\otimes \chi)$ simultaneously in the level, conductor and $s$ aspects. Assume that the form $f$ is of level $P$. Let $\vartheta=\log P/\log M$, and denote by $\mathcal{Q}_{f,\chi}=P M^2$ the (analytic) conductor of $L(1/2, f\otimes \chi)$. Blomer and Harcos's result implies 
\[L(1/2, f\otimes \chi)\ll \mathcal{Q}_{f,\chi}^{1/{4}+\varepsilon} \Big(  \mathcal{Q}_{f,\chi}^{-1/{(16+8\vartheta)}}+ \mathcal{Q}_{f,\chi}^{-{(1-\vartheta)}/{(8+4\vartheta)}}\Big )\] for any $0<\vartheta<1$. More recently, by averaging the fourth moment of $L( 1/2 , f \otimes \chi)$ over a family of newforms $f$, Khan \cite{Kh2} deduced that \[L(1/2,f\otimes \chi)\ll  \mathcal{Q}_{f, \chi}^{1/{4}+\varepsilon}\Big (\mathcal{Q}^{-{1}/{(8+4\vartheta)}}_{f,\chi}+ \mathcal{Q}^{-{(1-2\varpi)\vartheta}/{(16+8\vartheta  )}}_{f,\chi}\Big   ),\] where $\varpi$ denotes current best exponent towards the $
\GL(2)$ Ramanujan Conjecture. This, however, completely solved the hybrid subconvexity problem for ${\GL}(2)\times {\GL}(1)$ in the prime level and conductor aspects after	roughly thirteen years. 

In another direction, regarding the $L$-function $L(1/2,f\otimes g \otimes \chi)$, Sun \cite{Sun} proved a subconvex bound in the depth-aspect that
\[          L(1/2, g\otimes h \otimes \chi)\ll   {M^\prime}^{1-{1}/{16}+\varepsilon}    ,\]     where $g,h$ are any fixed holomorphic newforms, and the primitive character $\chi $ is of modulus $M^\prime$, $M^\prime=P^r$ for any integer $r\ge 12$. If the modulus is a prime $P$, quite recently, Ghosh \cite{Gh} attained 
\[          L(1/2, g\otimes h \otimes \chi)\ll  {P}^{1-1/23+\varepsilon}    .\]  
In addition, Sharma and Sawin \cite{Sh} were able to generalize to the higher rank case, by showing that, for any fixed {${\GL}(3)$} Hecek-Maa\ss\ form $\pi$ and {${\GL}(2)$} newform $g$, 
\[     L(1/2, \pi \otimes g\otimes \chi)\ll  P^{{3}/{2}-1/{16}+\varepsilon}    ,\] where the primitive character $\chi $ is still of prime modulus $P$. It is noticeable essentially they were working on Hecke cuspidal forms on {${\GL}(2)$} of level $P^2$ with the nebentypus $\chi^2$. Just lately, Kumar, Munshi and Singh \cite{KMS} investigated the hybrid subconvexity problem for the ${\GL}(3)\times {\GL}(2)$ $L$-functions in the levels aspects. Assume that the normalized form $\pi$ is of level $P$ and $g$ is of level $M$, respectively. They achieved that 
  \[L(1/2,\pi \otimes g)\ll    \mathcal{Q}_{\pi, g}^{1/{4}+\varepsilon}\Big (\mathcal{Q}_{\pi, g}^{-{(3-2\nu)}/{(16+24\nu)}}+ \mathcal{Q}_{\pi,g}^{-{(2\nu-1)}/{(16+24){\nu})}}\Big ) \] for any $1/2<\nu<3/2$, where $\mathcal{Q}_{\pi, g}=P^2 M^3 $ is the (analytic) conductor of $L(1/2,\pi \otimes g)$, and the exponent $\nu=\log P/\log M$. The
level aspect subconvexity problem for any genuine $ \GL(d)$ $L$-function with $d\ge 4$ still remains a very important
but wide-open problem.

Let $\Pi$ be a normalized Hecke-Maa\ss\ form on $\GL(4)$ of level $P$ (see \S2 for definition and backgrounds). In this paper, we shall consider subconvexity for the ${\GL}(4)\times {\GL}(1)$ $L$-function $L(s, \Pi\otimes \chi)$ in the level and conductor aspects simultaneously. To be precise, we are able to establish the following results:
\begin{theo} ~\label{58}
Let $P,M\ge 1$ be two primes such that $(P,M)=1$. Let $\chi$ be a primitive Dirichlet character modulo $M$. Let $\theta={\log P}/{\log M}$. Then, for any normalized Hecke-Maa\ss\ form $\Pi$ of level $P$ with trivial nebentypus, we have
\bee\label{2017102000101}\begin{split} &L (1/2, \Pi \otimes \chi) \ll \mathcal{Q}^{1/{4}+\varepsilon} \left  ( 
\mathcal{Q}^{-{(2- 5\theta)}/{3(32+8\theta ) }}+\mathcal{Q}^{-{\theta }/{(16+4\theta)  }} 
+\mathcal{Q}^{-{(1-2\mu)\theta }/{(8+2\theta)  }}  \right.\\
&\phantom{=\;\;}\left.   \hskip 9.36cm     
 +\mathcal{Q}^{(1-(4+2\mu))\theta/{3(16+4\theta)  }} \right )\end{split}
\ene for any $0<\mu<1/2$, where $\mathcal{Q}=PM^4$ denotes the (analytic) conductor of $L (1/2, \Pi \otimes \chi) $, and the implied $\ll$-constant depends only on $\varepsilon$ and the Langlands parameters $\alpha_i$, $1\le i\le 4$, as shown in \eqref{de040400444}.
\end{theo}
We thus obtain the hybrid subconvex bound in both parameters $P,M$ with $1/(4+2\mu)<\theta<2/5$ for any $0<\mu<1/2$. The main theorem above follows from an immediate application of an estimate upon the cancellations for the average of the coefficients of the Dirichlet series for $L(1/2,\Pi \otimes \chi)$.
\begin{theo} \label{03}
Let the notations be as in Theorem \ref{58}. Let $V$ be a smooth function, compactly supported on $[ 1/2 , 5/2 ]$ with bounded derivatives. Denote by $A_\Pi(n,1,1)$ the $n$-th coefficient of the Dirichlet series for $L(s,\Pi \otimes \chi)$. Then, for any $M^{1/(4+2\mu)}<P< M^{2/5}$ and $\mathcal{Q}^{1/2-\delta}\leq X\leq \mathcal{Q}^{1/2+\varepsilon}$ with $0<\mu <1/2$ and $\varepsilon,\delta>0$  being arbitrary, we have
\begin{align} \label{107}
\sum_{n\ge 1} \f{A_\Pi(n,1,1)\chi(n)}{\sqrt{n}}V\lf(\f{n}{X}\ri)\ll  X^{\fone{2}+\varepsilon} \lf (
\f{ \mathcal{Q}^{\delta} P^{5/8} }{M^{1/4}}+\fone{P^{1/2-\mu}}+\f{ \mathcal{Q}^{\delta} M^{1/4}}{P^{1+\mu/2}}+\fone{P^{1/4}} \ri),
\end{align}where the implied $\ll$-constant depends only on $\varepsilon$ and the Langlands parameters $\alpha_i$, $1\le i\le 4$.
\end{theo}
Indeed, observing that 
\[ \begin{split}L (1/2, \Pi \otimes \chi)&\ll \mathcal{Q}^\varepsilon \lf (\mathcal{Q}^{\fone{4}-\f{\delta}{2}}+ \sup_{  \mathcal{Q}^{1/2-\delta}\leq X\leq \mathcal{Q}^{1/2+\varepsilon}} \lf|\sum_{n\ge 1} \f{A_\Pi(n,1,1)\chi(n)}{\sqrt{n}}V\lf(\f{n}{X}\ri)\ri|\ri)\\
&\ll \mathcal{Q}^{\fone{4}+\varepsilon} \lf (\mathcal{Q}^{-\f{\delta}{2}} + \f{ \mathcal{Q}^{\delta} P^{5/8} }{M^{1/4}}+\f{ \mathcal{Q}^{\delta} M^{1/4}}{P^{1+\mu/2}}+\fone{P^{1/2-\mu}}+\fone{P^{1/4}}  \ri ) ,\end{split}\]one deduces \eqref{2017102000101}, upon equating the first three terms in the parentheses above. 

As a direct application of Theorem \ref{58} above, we obtain
\begin{coro}With the notations being as before, there exists a family of Hecke-Maa\ss\ forms $\Pi$ of level $P$, with $P\asymp M^{2/7+\varepsilon}$, such that
\[ L (1/2, \Pi \otimes \chi) \ll _{\Pi, \varepsilon} \mathcal{Q}^{1/{4}-1/60+\varepsilon} .\]\end{coro}

The basic manoeuvre in the paper follows a well-trodden but powerful path in studying $L$-functions involving levels aspects, namely the moment method plus the amplification. We shall employ a decomposition formula which was rooted in Munshi's work \cite{Mun21}, and then explicitly revealed by Holowinsky and Nelson \cite{HN}. This makes it feasible to attack the hybrid subconvexity problem for the $\GL(4)\times \GL(1)$ $L$-function by using the decomposition formula twice; otherwise, it seems impossible to succeed. It is crucial that at the first stage, we artificially introduce an amplifier and then we resort to this trick again before the second decomposition of the dualized $L$-values. To increase the degree of freedom of the $L$-structures, we gradually introduce the flexible parameters $\T$, $\R$, $\R^\ast$ and $\mathcal{S}$ relative to the fixed parameters $P,M$ attached to the central value $L(1/2,\Pi\otimes \chi)$. It, on the other hand, is very challenging to devise them to draw down on each other, so as to extract an admissible bound to beat the convexity barrier. 

To expose the method in this paper as clearly as possible, we have handled the degenerate term $\mathcal{S}^{\text{Deg.}}_{\Pi}$ in a relatively direct manner in \S3.3; this enables us lose the coverage of subconvexity for $\varepsilon<\theta<1/5$ (observing that morally the degenerate term $\mathcal{S}^{\text{Deg.}}_{\Pi}$ is always dominated by the non-degenerate term $\mathcal{S}^{^{\text{Non-de.}}}_{\Pi}$). Indeed, if one proceeds by repeating the argument as in \S3.1 in full details one can establish the range $\varepsilon<\theta<2/5$, where the Vorono\u{\i} summation formula in ramified case would be put into use instead. One might also wander if the method implies a way out of the subconvexity problem for $\GL(5)$. We are now working on this, and expect to report the progress in the future article.

Throughout the paper, $\varepsilon$ always denotes an arbitrarily small positive constant which might not be the same at each occurrence. 

\section{Prerequisites}
\subsection{Automorphic $L$-functions and functional equations}  In this part, we shall give a recap of the theory of automorphic $L$-functions. Let $N\ge 1$ be a square-free integer. We introduce the congruence subgroup \[\Gamma_1(N):=\text{Stab}_{{\rm{SL}(4, \Z)}}((0,0,0,1)) \subset {\rm{SL}}(4, \Z).\]
Let $\Pi$ be a normalized Hecke-Maa\ss\ form of type $v=(v_1,v_2,v_3)\in \C$ for the congruent subgroup $\Gamma_1(N)$ with trivial nebentypus, which has a Fourier-Whittaker expansion with the Fourier
coefficients $A_\Pi(m_1,m_2,m_3)$; see, e.g., \cite[Chapter 9]{Go} for definition and backgrounds. The Fourier coefficients of $\Pi$ and that of its contragredient $\widetilde{\Pi}$ are
related by $A_\Pi (m_1,m_2,m_3) = A_{\widetilde{\Pi}} (m_3,m_2,m_1)$ for any $(m_1,m_2,m_3, N) = 1$, with $A_\Pi (1, 1,1) = 1$. The Jacquet-Langlands $L$-function is given by $L(s,\Pi)=\sum_{n\ge 1}A_\Pi(n,1,1)n^{-s}$ for $\text{Re}(s)>1$. Now, let \bee \label{de040400444} \begin{split}&\alpha_1={3}/{2} -v_1-2v_2-3v_3,\quad 
\alpha_2=-{3}/{2} +3v_1+2v_2+v_3,\\
&\alpha_3=-{1}/{2} -v_1+2v_2+v_3,\quad
 \alpha_4={1}/{2} -v_1-2v_2+v_3\end{split}\ene be the Langlands parameters of $\Pi$. We define the completed $L$-function \[\Lambda(s, \Pi)=\mathcal{Q}^{s/2}_\Pi \varepsilon(\Pi) L_\infty (s,\Pi)L(s,\Pi),\]where 
\[L_\infty(s,\Pi)=\prod_{i=1}^4\Gamma_{\mathbb{R}}(s+\alpha_i),\] the (analytic) conductor $\mathcal{Q}_\Pi\asymp N$, and $|\varepsilon(\Pi)|=1$.  Here, we have followed the notational
convention that $\Gamma_{\mathbb{R}}(s)=\pi^{-s/2}\Gamma(s/2)$. One thus has the following functional equation \[\Lambda(s,\Pi)=\varepsilon(\Pi) \Lambda(1-s,
\widetilde{\Pi}).\] 

Next, let $M\ge 1$ be a square-free integer satisfying that $(N,M)=1$ by a little of abuse of notation. For any primitive Dirichlet character modulo $M$, we turn to considering the $\text{GL}(4)\times\text{GL}(1) $ twisted $L
$-function $L(s,\Pi\otimes \chi)$, which is given by the Dirichlet series $\sum_{n\ge 1}A_\Pi(n,1,1)\chi(n)n^{-s}$ for $\text{Re}(s)>1$. Now, one defines the Gamma factor by \[\begin{split}L_{\infty,\varrho} (s, \Pi\otimes \chi)&=\prod_{i=1}^4\Gamma_{\mathbb{R}}(s+\alpha_i+\varrho),
\end{split}\]where $\varrho=0$ or $1$ according to whether $\chi$ is even or odd. The completed $L$-function defined by \[\Lambda(s, \Pi\otimes \chi)=\mathcal{Q}^{s/2}  L_{\infty,\varrho} (s, \Pi\otimes \chi) L(s, \Pi\otimes \chi)\] with $\mathcal{Q}\asymp NM^4$ is thus an entire function with an analytic continuation to all $s\in \C$, and satisfies the function equation that \[\Lambda(s,\Pi\otimes \chi)=\varepsilon(\Pi\otimes \chi)\Lambda(1-s, \widetilde{\Pi}\otimes  \overline{\chi}).\]
\subsection{Vorono\u{\i} summation formula}Before stating the exact formula, we need define a hyper-Kloosterman sum of a special type which has appeared in the $\GL(4)$ Vorono\u{\i} formula. Let $a,c\in \Z$, and let
\[\mathbf{q}=(q_1,q_2), \quad \mathbf{d}=(d_1,d_2)\]
 be 2-tuples of positive integers satisfying the divisibility conditions
\[d_1|q_1c, \quad  d_2|\f{q_1q_2c}{d_1}.\]
The hyper-Kloosterman sum $\mathcal{KL}_2(n, m,c;\mathbf{q},  \mathbf{d})$ is defined as
\bee \label{342kl2}\mathcal{KL}_2(n, m,c;\mathbf{q},  \mathbf{d})=\smu_{x_1 \bmod \f{q_1c}{d_1} } \smu_{x_2 \bmod \f{q_1q_2c}{d_1d_2} } e{\lf (\f{d_1 x_1 n}{c}+\f{d_2\overline{x_1x_2}}{\f{q_1c}{d_1}}+\f{{x_2}m}{ \f{q_1q_2c}{d_1d_2}}\ri )} .\ene
 Let now $\omega$ be a smooth function compactly supported $[1/2,5/2]$ with bounded derivatives, and $\widetilde{\omega}$ its Mellin transform. For any given Maa\ss\ form $\Pi$ of prime level $P$, let the Langlands parameters be as in \eqref{de040400444}. For all $s\in \C$, define \[G_+(s)=\varepsilon(\Pi)  \prod_{i=1}^4\f{\Gamma_{\mathbb{R}}(1-s-\alpha_i)}{\Gamma_{\mathbb{R}}(s+\alpha_i)}, \quad  G_-(s)=\varepsilon(\Pi)  \prod_{i=1}^4\f{\Gamma_{\mathbb{R}}(2-s-\alpha_i)}{\Gamma_{\mathbb{R}}(s+\alpha_i+1)},\]
and set
\[\mho_{\pm} (x;\omega)=\fone{2\pi i} \int_{-\sigma} \widetilde{\omega}(s)x^s G_{\pm}(s)\ud s\] for $x>0$ and some $\sigma>0$. We then proceed by writing
\[\mho(x;\omega)=\mho_+(x;\omega)+\mho_{-}(x;\omega),\quad \mho(-x;\omega)=\mho_+(x;\omega)-\mho_{-}(x;\omega).\]
The $\GL(4)$ Vorono\u{\i} summation formula is thus given by
\begin{lem}\label{9945fyhjt3464lljyrdfv432}
Let $a\in \Z, c\in \mathbb{N}^+$ with $(c,aP)=1$. With notations as above, we have
\begin{align*}&\sum_{n\neq 0}  A_{\Pi}(n,1,1) \, 
 e{ \lf( \f{an}{c} \ri)  }
 \omega \lf (\f{n}{X} \ri)=c\sqrt{P}\sum_{\pm}\sum_{d_1|c}\sum_{d_2|\f{c}{d_1}}\sum_{m\neq 0} \f{\overline{A_\Pi(m, d_2,d_1)}}{|m|d_2d_1}\\
 &\hskip 7cm \mathcal{KL}_2(-\overline{aP}, m,c;(1,1) , (d_1,d_2) )\, \mho\lf ( \f{mXd^2_2 d^3_1}{c^4P}  ;\omega\ri)
.\end{align*}
\end{lem}
The assertion for Hecek-Maa\ss\ forms with trivial levels has already given in many works such as \cite{MS3,IT,CL}. The version involving the level aspect can be referred to Zhou \cite{Zh} who formulated the analog of the Vorono\u{\i} summation formula in the $\GL(3)$-case by appealing to the construction of a double Dirichlet series and
the functional equations of $L$-functions twisted by Dirichlet characters. Indeed, Corbett's formula \cite[Theorem 1.1]{C} can cover this, which proceeds in the
framework of the adelic interpretations; we will give the proof in the Appendix A.1.

We might proceed by writing $\mathcal{KL}_2(\overline{a}, m,c;(1,1) , (d_1,d_2) )$ asymptotically as
\[ \f{1}{d^2_1d_2}\smu_{\substack{x,y,z\bmod  c\\xyz \equiv 1\bmod c}}\e{\f{d_1 x \overline{a}+ d_1d_2 m {y  }+ d_1d_2z }{c}}.\]With this, one has seen the relation that
\bee \label{9xj3c43c4xe}
\begin{split}&\sum_{n \neq 0}  \f{A_{\Pi}(n,1,1)}{\sqrt{n}} \, 
 e{ \lf( \f{an}{c} \ri)  }
 \omega \lf (\f{n}{X} \ri)\\
 &\hskip2cm \Longrightarrow \,  \sum_{d_1|c}\sum_{d_2|\f{c}{d_1}}\sum_{m \neq 0} \f{\overline{A_\Pi(m, d_2,d_1)}}{\sqrt{m d^2_2d^3_1}}\, \f{\Kl(-\overline{aP} ,m;d_1,d_2,{c})} {c}  \, \mho \lf ( \f{m Xd^2_2 d^3_1}{c^4P}  ;\omega\ri)
.\end{split}
\ene
 where the sum \bee \label{33xc2v4c2}{\Kl(m ,n;d_1,d_2,{c})} =\smu_{\substack{x,y,z\bmod  c\\ xyz\equiv 1 \bmod {c}}}\e{ \f{d_1 m x + d_1d_2 n {y  }+ d_1d_2z}{c}}\ene for any integers $m,n,c\ge 1$.
 
We have a good control of the asymptotic behavior of $\Omega_\pm$, upon expressing it in term of the oscillatory integral; see, e.g., \cite[Lemma 5.2]{CL}. 
\begin{lem}\label{3453464lljyrdfv432} For any fixed integer $\mathcal{K}\ge 1$ and $x\gg 1$, we have
\begin{align*}&\mho_+(x;\omega)=x\int_{\mathbb{R}^+} \omega(y) \sum_{j=1}^{\mathcal{K} } \fone{(xy)^{{j}/{4}+1/{8}} }\Big \{ c_j  e{\lf (4(xy)^{1/{4}} \Big  )} \mathcal{W}^+_j   \Big (8\pi (xy)^{{1}/{4}} \Big )
  \\
  &  \hskip 5cm
 +d_j  e{\Big  (-4(xy)^{1/{4}} \ri )} {\mathcal{W}^-_j  } \Big (8\pi (xy)^{{1}/{4}} \Big )   \Big  \}\ud y+O\Big( x^{{(3-\mathcal{K})}/{4}-1/{8}}\Big ),\end{align*} 
 where $\mathcal{W}^\pm_j $ are certain  smooth weight functions satisfying that $x^i{\mathcal{W}^\pm_j }^{(i)}(x)\ll 1$ for any $i\ge 1$, and $c_j,d_j$ are suitable constants depending on the four parameters $\alpha_i$, $1\le i\le 4$, above. Furthermore, $\mho_-(x;\omega)$ has the same expression except values of constants.
\end{lem}

\section{ Proof of Theorem \ref{03}} 
In this section, we are dedicated to proving Theorem \ref{03}. The point of departure, however, is the assembly with the amplification method. For any fixed $ T\ge  1$, we introduce a set of moduli $\mathfrak{t}$:
\[\begin{split} &\mathfrak{t}=\{T\le t\le 2T  \text{ is prime, } (t,MP)=1 \}.\end{split}\]
We will choose the sequence $\{\gamma _t\}_t \in \C$ of moduli at most one, supported on the set $\mathfrak{t}$, such that \[\fone{\T}\sum_{t\in \mathfrak{t} } \gamma_t \overline{\chi}(t)=1 .\]  Here, $ \T$ stands for the cardinalities of the set $\mathfrak{t}$, whcih satisfies that $\T \asymp T/\log T$. Notice that the sequence $\gamma _t$ is well-defined. For example, one can take $\gamma_t=\chi(t)$, whenever $t\in \mathfrak{t}$, and zero otherwise. Now, if one defines \bee \label{3432423512514} \mathcal{S}_{\Pi}(X,M,P)=\sum_{n\ge 1 }\f{A_\Pi(n,1,1)\chi(n)}{\sqrt{n}} V\lf (\f{n}{X} \ri)\ene with $V$ being a smooth function supported on $[1/2, 5/2]$ of bounded derivatives, we thus transform $ \mathcal{S}_{\Pi}(X,M,P)$ into 
\[ \mathcal{S}_{\Pi}(X,M,P)=\fone{ \T} \sum_{t\in \TT}  \gamma_t \sum_{n\ge 1} \f{A_\Pi(n,1,1)\chi (\overline{t} n) }{\sqrt{n}} \V{\f{n}{X}}. \] 
For any given primitive Dirichlet character $\chi \bmod M$, let us recall that Holowinsky and Nelson \cite{HN} has formulated a decomposition for $\chi$\begin{align} \label{20170389229} \chi(n)=\f{M}{\R\tau(\overline{\chi})} \sum_{r\ge 1} \chi(r)\e{\f{n\overline{r}}{M}} \U{\f{r}{\R}}-\fone{\tau(\overline{\chi})} \sum_{r 
\ge 1} {\rm{Kl}}_{2,\chi} (r,n;M)\, \widehat{U}\lf (\f{r}{M/\R}\ri),\end{align}
where $0<\R\le M$ is a freely chosen parameter. Here, the twisted Kloosterman sum 
\[{\rm{Kl}}_{2,\chi} (m,n;q)=\smu_{\substack{
x,y\bmod q\\xy \equiv 1 \bmod q}} \chi(x)e{\lf ( \f{mx+ny}{q}\ri)}\]for any integers $m,n,q\ge1$, and $ \widehat{U}$ is the Fourier transform of the smooth weight function $U$, i.e., $\widehat{U}(x)=\int_{\mathbb{R}^+}U(y)e(xy)\ud y$, which is supported on $[1/2,5/2]$ and normalized such that $ \widehat{U}(0)=1$. Recall that typical transformation formula for $\chi$ should be 
 \[\chi(n)=\f{1}{\tau(\overline{\chi})} \smu_{r \bmod M} \chi(r)\e{\f{n\overline{r}}{M}}.
 \]
One might 
recognize the first term on the right-hand side of \eqref{20170389229} as  a truncated approximation for $\chi$ by an exponential sum of the length less than the conductor; while, the second term should be thought of as its dual form. Now, an application of \eqref{20170389229} splits $ \mathcal{S}_{\Pi}(X,M,P)$ into three parts: \bee \label{jdid33d33324} \mathcal{S}_{\Pi}(X,M,P)= \mathcal{S}^{\text{Non-de.}}_{\Pi}(X,M,P)- \mathcal{S}^{\text{Dua.}}_{\Pi}(X,M,P)+ \mathcal{S}^{\text{Deg.}}_{\Pi}(X,M,P),\ene where \begin{align}
 \label{202010}\begin{split}
 \mathcal{S}^{\text{Non-de.}}_{\Pi}(X,M,P)=&\f{M}{ \mathcal{T} \R\, \tau(\overline{\chi})} \sum_{t\in \TT} \gamma_t \sum_{\substack{r\in \Z_{\neq 0} \\
(r,P)=1} }\chi(r)\U{\f{r}{R}} \\
&\hskip2.6cm \sum_{n\in \Z_{\neq 0} } \f{A_\Pi(n,1,1)}{\sqrt{n}}  \e{\f{ns\overline{tr}}{M}}  \V{\f{n}{X}},\end{split}
\end{align} 
\begin{align}\label{202020} \begin{split}& \mathcal{S}^{\text{Dua.}}_{\Pi}(X,M,P)=\f{1}{ \T\, \tau(\overline{\chi})}\sum_{\substack{t\in \TT}}\gamma_t \sum_{r\ge 1}\widehat{U}\lf (\f{r}{M/\R}\ri)\\
&\hskip4.6cm  \,\sum_{n\ge 1 } \f{A_\Pi(n,1,1)}{\sqrt{n}}\,
 {\rm{Kl}}_{2,\chi}(r,\overline{t}n;M)  \, \V{\f{n}{X}},\end{split}\end{align}
and  
\begin{align} \label{202010456247868}\begin{split}
 \mathcal{S}^{\text{Deg.}}_{\Pi}(X,M,P)=&\f{M}{ \T\R\,\tau(\overline{\chi})} \sum_{t\in \TT} \gamma_t \sum_{\substack{r\ge 1\\
P|r} }\chi(r)\sum_{n\ge 1} \f{A_\Pi(n,1,1)}{\sqrt{n}} \e{\f{n\overline{tr}}{M}} \U{\f{r}{\R}} \V{\f{n}{X}},\end{split}\end{align} 
respectively. Observe that $ \mathcal{S}^{\text{Deg.}}_{\Pi}$ survives, if and only if $\R > P$. We will assume henceforth that the parameters $\R$ and $\T$ satisfy that
\bee \label{023e342}
 \R<\T<M.
\ene  
\subsection{Treatment of $ \mathcal{S}^{\text{Non-de.}}_{\Pi}$}
This subsection is devoted to estimating $ \mathcal{S}^{\text{Non-de.}}_{\Pi}$. Recall \eqref{202010}. By the reciprocity law, one writes \[\e{   \f{n\overline{tr}}{M}  }=\e{  - \f{n\overline{M}}{tr}  }\e{   \f{n}{trM}  } .\] 
Now, if one chooses  \bee \label{234235451}  \R=\f{X}{M\T}, \ene the harmonic $e(n/(trM))$ is automatically {flat}. We thus turn to the following 
\bee \label{234233456341}\f{\sqrt{M}}{\T\R} \sum_{\substack{t\in \mathfrak{t}}}  \gamma_t \sum_{\substack{r\ge 1 \\(r,P)=1}} \chi(r) \sum_{n\ge 1 } \f{\A}{\sqrt{n}} \e{  - \f{n\overline{M}}{tr}  } \U{\f{r}{\R}} \V{\f{n}{X}}\ene in estimating the non-degenerate term $ \mathcal{S}^{\text{Non-de.}}_{\Pi}$. Our strategy next is to apply the $\rm{GL}(4)$-Vorono\u{\i} 
formula to the $n$-sum above, followed by the delicate analysis of the resulting structures
reformulations. This (essentially) produces a transformation for the display in \eqref{234233456341} as follows
\begin{align}
 \label{4040404000456456451} & \f{\sqrt{M}}{(\T\R)^2}\sum_{t\in \TT}  \gamma_t \sum_{\substack{ r \ge 1\\ (r,P)=1}} \chi(r )\sum_{d_1|tr}\,\sum_{d_2|tr/d_1}\,\sum_{m\ge 1} \f{\overline{A_\Pi(m,d_2,d_1)}}{\sqrt{md^2_2d^{3}_1}}  
 \notag \\
 &\hskip 5.6cm 
{\Kl(\overline{P}M, n;d_1,d_2,{tr})}\, {\mathscr{W}}_\pm\lf ( \f{m  d^2_2d^{3}_1}{\mathscr{N}},\f{r}{\R}\ri).
\end{align} Here, for any $r,\iota\in \mathbb{R}^+$, the smooth weight ${\mathscr{W}}_\pm$ is defined as ${\mathscr{W}}_\pm{ ({\iota},r)}= \mho_{\pm}{( \iota  ;{V})}  U(r
) $, and $\mathscr{N}=(\T\R)^4P/X^{1-\varepsilon}$. Upon combining with Lemma \ref{3453464lljyrdfv432}, one sees that the resulting integral transform $ \mho_{\pm}(x;V)$ asymptotically equals \bee\label{094v42}
 x ^{{5}/{8}}    \sum_{\pm}\int_{\mathbb{R}^+} \f{  V(y)  }{  y^{{3}/{8} }}
 e{\Big (\pm    4 ( xy ) ^{1/{4}  }  \Big)}
 \,
   {\mathcal{W}^\pm_1 } \Big (    8\pi ( xy  ) ^{1/{4}  }   \Big ) \ud y \ene plus an acceptable error term, with the weights $\mathcal{W}^\pm_1  $ satisfying that $ x^j {\mathcal{W}^\pm_1}^{(j)} (x)\ll X^\varepsilon$ for any $j\ge 1$. Recall the definition of ${\Kl(m, n;d_1,d_2,{tr})}$ as shown in \eqref{33xc2v4c2}. It suffices to take care of the salient case where $d_2=1$ in the following, the portion with $d_2$ being large contributes less significantly (as one expects) by a similar analysis. Note that $(t,r)=1$ by the pre-condition \eqref{023e342}; this implies that there are three options for $d_1$, that is, $d_1=1$, $d_1=t$ and $d_1|r,d_1>1$ ($d_1=tr$ cannot occur, in view of that $\mathscr{N}<(\T\R)^3X^{-\varepsilon}$). For clarity of presentation, it suffices to deal with the case of $d_1=1$ (this will be responsible for the dominated contribution of $ \mathcal{S}^{\text{Non-de.}}_{\Pi}$); the same argument works for other situations which brings relatively small contributions. In this sense, the task is essentially boiled down to 
\begin{align}
 \label{4040404000456456451}  \f{\sqrt{M}}{(\T\R)^2}\sum_{t\in \TT}  \gamma_t \sum_{\substack{ r \ge 1\\ (r,P)=1}} \chi(r )\,\sum_{n\ge 1} \f{\overline{A_\Pi(n,1,1)}}{\sqrt{m}}  
 \,
 {\Kl(\overline{P}M, n;{tr})} \,{\mathscr{W}}_\pm\lf ( \f{n  }{\mathscr{N}},\f{r}{\R}\ri),
\end{align} where the classic hyper-Kloosterman sum $\Kl$ is defined as \[\Kl(m, n;{q})=\smu_{\substack{x,y,z\bmod  q\\ xyz\equiv 1 \bmod {q}}}\e{ \f{m x +  n {y  }+z}{q}}\] for any integers  $m,n,q\ge 1$. Observe that $(M,r)=1$, on account of the assumption that $\R<M$. Let $n=n_0n_1$ with $n_0|(tr)^\infty$ and $(n_1,tr)=1$. By a standard computation (see, e.g., \cite[Section 2]{Sm}), one might write \[ {\Kl(\overline{P}M, n;{tr})} =\fone{\varphi(tr)}\sum_{\psi \bmod {tr}} \tau^\star(\overline{\psi},n_0) \tau^2(\overline{\psi}) \psi (\overline{P}M)\psi(n),\]
with \[\tau^\star(\overline{\psi},n_0)=\sum_{\alpha\bmod {tr}} e{\lf (\f{n_0 \alpha}{tr} \ri )} \overline{\psi}(n_0 
\alpha).\] With the help of this, the expression \eqref{4040404000456456451} becomes
 \begin{align}
 \label{4040404000456456451}  &\f{\sqrt{M}}{(\T\R)^2}\sum_{t\in \TT}\gamma_t \sum_{\substack{ r \ge 1\\ (r,P)=1}} \chi(r) \sum_{\psi \bmod {tr}}\f{\tau^2(\overline{\psi})  \psi (\overline{P}M)}{ \varphi(tr)}   \notag\\
 &\hskip2.7cm\sum_{\substack{n\ge 1\\ n=n_0n_1\\n_0|(tr)^\infty\\(n_1,tr)=1}}   \tau^\star(\overline{\psi},n_0) 
\f{\overline{A_\Pi(n,1,1)} \psi(n)}{\sqrt{n}}  {\mathscr{W}}_\pm\lf ( \f{n  }{\mathscr{N}},\f{r}{\R}\ri).
\end{align}  
Now, we proceed by introducing a new amplifier $\beta_s $. For any $ S\ge  1$, we define another set of moduli $\mathfrak{s}$ by
\[\begin{split} &\mathfrak{s}=\{{S}\le s\le 2{S} \text{ is prime}  ,\, (s,MP)=1,  \text{ and } (s,t)=1\text{ for any } t\in \mathfrak{t}   \},\end{split}\] with $\sharp{\mathfrak{s}}=\mathcal{S}\asymp S/\log S$ to be chosen optimally in the sequel. Here, as before the sequence $\{\beta_s\}_s$ is of moduli at most one, supported on the set $\mathfrak{s}$ and satisfies that  \[\fone{\mathcal{S}}\sum_{s\in \mathfrak{s} } \beta_s \chi(s)=1   .\]  
By an assembly with this new amplifier, the multiple sum \eqref{4040404000456456451} is transformed into
 \begin{align*}
&\f{\sqrt{M}}{\mathcal{S}(\T\R)^2 }\sum_{s\in \ST,t\in \TT}\beta_s \gamma_t \sum_{\substack{ r \ge 1\\ (r,P)=1}} \chi(r) \sum_{\psi \bmod {tr}}\f{\tau^2(\overline{\psi})\psi (\overline{P}M)}{\varphi(tr)}\\
&\hskip 4cm  \sum_{\substack{n\ge 1\\ n=n_0n_1\\n_0|(tr)^\infty\\(n_1,tr)=1}}   \tau^\star(\overline{\psi},n_0) \f{\overline{A_\Pi(n,1,1)}  \psi(n\overline{s})}{\sqrt{n}} {\mathscr{W}}_\pm\lf ( \f{n  }{\mathscr{N}},\f{r}{\R}\ri)
\end{align*} 
The manoeuvre for the next stage is to resort to the formula \eqref{20170389229} twice. Let $1\le \R^\ast\le \T\R$ be a parameter which will be explicitly determined later. The display above thus turns into
  \begin{align*}
 \label{4040404000456456451}  &\f{\sqrt{M}}{\mathcal{S}(\T\R)^2 }\sum_{s\in \ST,t\in \TT}\beta_s \gamma_t \sum_{\substack{ r \ge 1\\ (r,P)=1}} \chi(r)  \sum_{\psi \bmod {tr}}\f{\tau^2(\overline{\psi}) \psi (\overline{P}M)}{\varphi(tr)}\,
\\
 & \hskip 1cm  \sum_{\substack{n\ge 1\\ n=n_0n_1\\n_0|(tr)^\infty\\(n_1,tr)=1}}   \tau^\star(\overline{\psi},n_0) \,
  \f{\overline{A_\Pi(n,1,1)}}{\sqrt{n}} \lf \{\f{tr}{\R^\ast \tau(\overline{\psi})}
 \sum_{r^\ast \ge 1} \psi(r^\ast)  e{\lf ( \f{n \overline{sr^\ast}}{tr}\ri)} {\mathscr{V}}\lf ( \f{r^\ast }{\R^\ast } \ri)\\
 &\hskip 5.2cm 
-\fone{\tau(\overline{\psi})} \sum_{\nu\ge 1 }     {\rm{Kl}}_{2, \psi} (\nu, n\overline{s} ;tr  ) \, \widecheck{\mathscr{V}}\lf ( \f{\nu \R^\ast }{\T\R} \ri)
    \ri \}{\mathscr{W}}_\pm\lf ( \f{n  }{\mathscr{N}},\f{r}{\R}\ri),
\end{align*} 
where $ \mathscr{V}$ is a Schwarz function, which behaves like $U$, and $\widecheck{\mathscr{V}}$ denotes its Fourier transform given by $\widecheck{\mathscr{V}}(x)=\int_{\mathbb{R}^+}{\mathscr{V}}(y)e(xy)\ud y$ such that $\widecheck{\mathscr{V}}(0)=1$. Observe that 
\[\begin{split}&\fone{\varphi(tr)} \sum_{\psi \bmod {tr}} \tau^\star(\overline{\psi},n_0) \tau(\overline{\psi}) \psi(r^\ast \overline{P}M)\\
&=\fone{\varphi(tr)} \smu_{x,y \bmod {tr}} e{\lf ( \f{n_0 x+y}{tr}  \ri)} \sum_{\psi \bmod {tr}} \overline{\psi}(n_0xy)\psi (r^\ast \overline{P}M)\\
&={\rm{Kl}}_2(r^\ast \overline{P}M,1;tr).\end{split}\]
We are thus led to the following two terms 
  \begin{align*}
& \mathcal{S}^{\text{Non-de.,1}}_{\Pi}(X,M,P)=\f{\sqrt{M}}{\mathcal{S} \T \R \R^\ast } \sum_{s\in \ST,t\in \TT}\beta_s \gamma_t \sum_{\substack{ r \ge 1\\ (r,P)=1}} \chi(r)  \sum_{r^\ast \ge 1 }  \sum_{n\ge 1 }   \f{\overline{A_\Pi(n,1,1)} }{\sqrt{n}} \\
&\hskip 5.2cm e{\lf ( \f{n \overline{sr^\ast}}{tr}\ri)}
{\rm{Kl}}_2(\overline{P}M r^\ast,1;tr)\,{\mathscr{V}}\lf ( \f{r^\ast }{\R^\ast } \ri){\mathscr{W}}_\pm\lf ( \f{n  }{\mathscr{N}},\f{r}{\R}\ri),
\end{align*} 
and 
  \begin{align}\label{904k30fkff}
& \mathcal{S}^{\text{Non-de.,2}}_{\Pi}(X,M,P)=\f{\sqrt{M}}{\mathcal{S}(\T\R)^2} \sum_{s\in \ST,t\in \TT}\beta_s \gamma_t \sum_{\substack{ r \ge 1\\ (r,P)=1}} \chi(r)  \sum_{\nu \ge 1}  \, \sum_{n\ge 1 }   \f{\overline{A_\Pi(n,1,1)} }{\sqrt{n}} \notag\\
&\hskip 7cm  \mathfrak{C}(s,\nu, n;tr)\,\widecheck{\mathscr{V}}\lf ( \f{\nu \R^\ast }{\T\R} \ri) {\mathscr{W}}_\pm\lf ( \f{n  }{\mathscr{N}},\f{r}{\R}\ri)
\end{align}  with 
\bee \label{0932jfddne2j2e2}\mathfrak{C}(s,m, n;q)=\smu_{\gamma \bmod {q}}  e{\lf ( \f{m{\gamma}+ n\overline{s\gamma}}{q}\ri)} {\rm{Kl}}_2(\overline{P}M \gamma,1 ;q)\ene for any integers $s,m,n,q\ge1$ with $(s,q)=1$.

In the rest of the paper, we will continue to assume that
\bee \label{023eddwd}
\R<\mathcal{S}, \quad \mathcal{S}+\T<\R^\ast<\T \R.
\ene
\subsubsection{Bounding $ \mathcal{S}^{\text{Non-de.,1}}_{\Pi}$} Let us now begin with estimating the term $ \mathcal{S}^{\text{Non-de.,1}}_{\Pi}$. By applying the $\rm{GL}(4)$-Vorono\u{\i} 
formula, together with Lemma \ref{3453464lljyrdfv432}, this term turns out to be bounded by 
\begin{align}
 \label{c2vovir} 
&\f{\sqrt{M}}{\mathcal{S}(\T\R)^2\R^\ast } \sum_{s\in \ST,t\in \TT}\beta_s \gamma_t \sum_{\substack{ r\ge 1\\ (r,P)=1}} \chi(r)  \sum_{r^\ast \ge 1}  \sum_{d_1|tr}\,\sum_{d_2|tr/d_1}\,\sum_{m \ge 1} \f{A_\Pi(m,d_2,d_1)}{\sqrt{md^2_2d^{3}_1}}  \, 
\notag\\
& \hskip2.6cm\Kl(-sr^\ast\overline{P}, m;d_1,d_2, tr) \,
 {\rm{Kl}}_2(\overline{P}M r^\ast , 1;tr)\,
{\mathscr{V}}\lf ( \f{r^\ast }{\R^\ast } \ri) \,\widetilde{{\mathscr{W}}_\pm}\lf ( \f{m d^2_2d^{3}_1 }{X},\f{r}{\R}\ri),
\end{align}
where 
\begin{align*}\widetilde{{\mathscr{W}}_\pm}(x,r)&=U(r)x^{5/8}\int_{1}^{10X^{100\varepsilon}}\int_{\mathbb{R}^+}\f{V(y)t^{1/4}}{y^{3/8}}\,e(4(ty)^{1/4}\pm 4(tx)^{1/4})\\
&\hskip 6cm \mathcal{W}_1(8\pi(ty)^{1/4} ) \,\mathcal{W}^\pm_1(8\pi(tx)^{1/4} ) \ud y\ud t\end{align*}
for any $x,r\in \mathbb{R}^+$. One might identify $\widetilde{{\mathscr{W}}_\pm}$ as a Schwarz function, which behaves like one with the partial derivatives $x^i r ^j
 \f{\partial ^i}{\partial x^i} 
  \f{\partial ^j} {\partial r^j} \widetilde{{\mathscr{W}}_\pm}(x,r)\ll X^\varepsilon $ for any $i,j\ge 0$. As hinted before, for the purpose of gaining an effective control over the magnitude of the sum \eqref{c2vovir}, the main focus is enough to be on the focal situation where $d_1=d_2=1$. The expression in \eqref{c2vovir} is dominated by
 \begin{align*} 
&\f{\sqrt{M}}{\mathcal{S}(\T\R)^2\R^\ast } \sum_{s\in \ST,t\in \TT}\beta_s \gamma_t \sum_{\substack{ r \ge 1\\ (r,P)=1}} \chi(r)  \sum_{r^\ast \ge 1}  \,\sum_{m \ge 1} \f{A_\Pi(m,1,1)}{\sqrt{m}}  \Kl(-sr^\ast\overline{P},m; tr)\, 
\\
& \hskip7cm
 {\rm{Kl}}_2(\overline{P}M r^\ast , 1;tr)\,
\,{\mathscr{V}}\lf ( \f{r^\ast }{\R^\ast } \ri) \widetilde{{\mathscr{W}}_\pm}\lf ( \f{m }{X},\f{r}{\R}\ri).
\end{align*}
We now introduce a smooth weight by letting $\phi(x)\in \mathcal{C}^{\infty}_0(\mathbb{R},\mathbb{R}^+)$ be a function supported on $[X^{-\varepsilon}/2, 5X^{\varepsilon}/2]$, identical to one for $x\in [X^{-\varepsilon}, 2X^{\varepsilon}]$, with the bounded derivatives. Via interchanging the order of sums and then applying the Cauchy-Schwarz inequality to get rid of Fourier coefficients $A_\Pi(m,1,1)$ with the help of the bound of ``Ramanujan on
average" $ \sum_{m\le X}|A_\Pi(m,b,c)|^2\ll P^\varepsilon (bcX)^{1+\varepsilon}$ (see, e.g., \cite[Section 3]{CL}), we arrive at 
\bee \label{nxji32w2212} \mathcal{S}^{\text{Non-de.,1}}_{\Pi}(X,M,P) \ll \f{\sqrt{M}}{(\T\R)^{2-\varepsilon}\mathcal{S}\R^\ast }   \, \mathcal{A}^{1/{2}}(X,P,M),\ene
where  $\mathcal{A}(X,P,M)$ is given by
 \begin{align*}  &\sum_{m\ge 1} \phi {\lf ( \f{m}{XL} \ri)}
\sum_{\substack{s_1,s_2\in \ST\\ t_1,t_2 \in \TT}}\beta_{s_1}\overline{\beta_{s_2}} \gamma_{t_1}  \overline{\gamma_{t_2}}\sum_{\substack{ r_1,r_2 \ge 1\\ ( r_1r_2 ,P)=1}} \chi( r_1)\overline{ \chi( r_2 ) }    \notag
 \\
& \hskip1.1cm \sum_{r^\ast_1,r^\ast_2 \ge 1} \Kl(-s_1r^\ast_1\overline{P}, m; t_1r_1)\,  {\rm{Kl}}_2(\overline{P}M r^\ast_1 , 1;t_1r_1)\,
\overline{\Kl(-s_2r^\ast_2 \overline{P}, m; t_2r_2)  }  \notag
\\
&\hskip 2.7cm  \overline{  {\rm{Kl}}_2(\overline{P}M r^\ast_2 , 1;t_2r_2)}\, 
{\mathscr{V}}\lf ( \f{r^\ast_1 }{\R^\ast } \ri) \widetilde{{\mathscr{W}}_\pm}\lf ( \f{m }{X},\f{r_1}{\R}\ri)\, \overline{
{\mathscr{V}}\lf ( \f{r^\ast _2}{\R^\ast } \ri) \widetilde{{\mathscr{W}}_\pm}\lf ( \f{m }{X},\f{r_2}{\R}\ri)}.
\end{align*}
Here, $L\ge 1$ is a parameter to be chosen after a while. This is reminiscent of the `Maass transfer trick' given by Munshi \cite{M27} and the trick of Blomer and Leung \cite[Section 5]{BL}. We are in a position to exploit Poisson summation with the modulus $t_1t_2r_1r_2$ to the sum over $n$. Note that $\R\ll \T$; this implies that $(t_1t_2, r_1r_2)=1$. The $m$-sum is essentially converted into
\begin{align*}&\f{X^{1+\varepsilon}L}{t_1t_2r_1r_2}\sum_{|h|\ll ( t_1t_2r_1r_2)^{1+\varepsilon}/XL }\,
 \sum_{\gamma \bmod{t_1t_2r_1r_2}} \Kl(-\overline{P} s_1r^\ast_1,\gamma ; t_1r_1)\\
 &\hskip 5.7cm  
\overline{\Kl(-\overline{P} s_2r^\ast_2, \gamma ; t_2r_2)  }\, e{\lf ( -\f{\gamma  h }{t_1t_2r_1r_2}\ri)}\end{align*}We will choose $L=P^{1/2+100\varepsilon}$, so that the non-zero frequences do not exist, upon recalling \eqref{234235451} which shows that $M^{1-\varepsilon }\sqrt{P}\ll \T \R\ll M^{1+\varepsilon }\sqrt{P}$. With this, it is necessary that $t_1=t_2$, and, moreover, $(r_1,r_2)=1$ which implies $r_1=r_2$; otherwise, the sum is (negligibly) small, on account of the fact that $\sum_{t \bmod q }\Kl(\ell t,1; q)$ essentially `vanishes' for any integers $\ell, q\ge 1$. Assume that $t_1=t_2=t$ and $r_1=r_2=r$, say. By Poisson summation with the modulus $tr$, it follows that the quantity we are faced with actually is the following
 \begin{align} \label{dej393jjs}
& \f{X^{1+\varepsilon}\sqrt{P}}{\T\R}\sum_{\substack{t \in \TT\\ s_1,s_2\in \ST}}|\beta_{s_1} \overline{\beta_{s_2}} \gamma_{t}  |^2\, \sum_{\substack{ r \ge 1\\ ( r ,P)=1}}  \sum_{r^\ast_1,r^\ast_2 \ge 1}   {\rm{Kl}}_2(\overline{P}M r^\ast _1, 1;tr) \, \overline{ {\rm{Kl}}_2(\overline{P}M r^\ast_2 , 1;tr)} \notag\\
& \hskip 4cm  \,
\sum_{\jmath \bmod {tr}}   
{\Kl(-\overline{P} s_1r^\ast_1, \jmath ; tr)  }\,
\overline{{\Kl(-\overline{P}s_2r^\ast_2, \jmath ; tr)  }} \, {\mathscr{V}}^\flat \lf ( \f{r^\ast_1 }{\R^\ast } , \f{r^\ast_2 }{\R^\ast } , \f{r}{\R } \ri)  ,
\end{align}
with $ {\mathscr{V}}^\flat(x,y,z)={\mathscr{V}}(x){\mathscr{V}}(y)U(z)$ for any $x,y,z\in \mathbb{R}^+$.

---\emph{When $s_1=s_2$ and $r^\ast_1=r^\ast_2$}. Note that $r^\ast_1 s_1\equiv r^\ast_2 s_2 \bmod {tr}$. If $r_1=r_2$, then $s_1=s_2$, upon recalling the pre-condition \eqref{023eddwd}. In this case, the sum $\mathcal{A}(X,P,M)$ is estimated as
 \bee \label{c42t25253155}\begin{split} 
& \ll \f{X^{1+\varepsilon}\sqrt{P}}{TR}\sum_{s\in \ST,t \in \TT}|\beta_{s}\gamma_{t}  |^2\sum_{\substack{ r \sim\R}}  \sum_{r^\ast \sim \R^\ast}  | {\rm{Kl}}_2(\overline{P}M r^\ast , 1;tr) |^2 \,    
\sum_{\jmath \bmod {tr}}   
|{\Kl(-\overline{P} sr^\ast, \jmath ; tr)  }|^2\\
&\ll  X^{1+\varepsilon}\sqrt{P} \mathcal{S}\R^\ast (\T\R)^4.
\end{split}\ene

---\emph{When $s_1\neq s_2$ and $r^\ast_1\neq r^\ast_2$}.  Next, we consider $r^\ast_1\neq r^\ast_2$; this implies $s_1\neq s_2$. The display \eqref{dej393jjs} is controlled by
 \[\begin{split} 
& \ll {X^{1+\varepsilon}\sqrt{P}}{(\T\R)^2}\sum_{s_1,s_2\in \ST,t \in \TT}|\beta_{s_1} \overline{\beta_{s_2}} \gamma_{t}  |^2\,\sum_{\substack{ r \ge 1}} \,   
 \sum_{r^\ast_1,r^\ast_2 \ge 1}   {\rm{Kl}}_2(\overline{P}M r^\ast _1, 1;tr) 
 \\
& \hskip4.1cm  \overline{ {\rm{Kl}}_2(\overline{P}M r^\ast_2 , 1;tr)}  \, \mathbf{1}_{r^\ast_1 \equiv r^\ast_2 s_2 \overline{s_1} \bmod {tr}}\, {\mathscr{V}}^\flat \lf ( \f{r^\ast_1 }{\R^\ast } , \f{r^\ast_2 }{\R^\ast } , \f{r}{\R } \ri) .
\end{split}\]
Observe that
\begin{align*}  & \sum_{r^\ast_2 \ge 1}   {\rm{Kl}}_2(\overline{P}M r^\ast_2 s_2\overline{s_1}, 1;tr)  \overline{ {\rm{Kl}}_2(\overline{P}M r^\ast_2 , 1;tr)}  \,\sum_{\substack{r^\ast_1\ge 1\\
r^\ast_1 \equiv r^\ast_2 s_2 \overline{s_1} \bmod {tr}}} \, {\mathscr{V}}\lf ( \f{r^\ast _1}{\R^\ast }  \ri)\,\, {\mathscr{V}}\lf ( \f{r^\ast _2}{\R^\ast }  \ri)\\
&\ll  \f{\R^\ast}{(tr)^{1-\varepsilon}} \sum _{|\ell |\ll {\T\R/{\R^\ast}^{1-\varepsilon}}}
\, \lf |
 \sum_{\hbar \bmod {tr}}  {\rm{Kl}}_2(\overline{P}M \hbar  s_2\overline{s_1}, 1;tr) \, \overline{ {\rm{Kl}}_2(\overline{P}M \hbar  , 1;tr)}\, e{\lf (\f{\ell \hbar }{tr}\ri)} \ri|
\end{align*}by invoking Poisson summation with the modulus $tr$ to the $r^\ast_2$-sum. Now, if $\ell=0$, this shows that $s_1\equiv  s_2 \bmod {tr}$ and thence $s_1=s_2$. We thus return to the former situation. We will be devoted to the non-zero frequences in what follows. By the estimate for the triple sum \bee \label{x431cr34pov2}\sum_{\gamma \bmod {q}}  {\rm{Kl}}_2(m \gamma , 1;tr)\, { {\rm{Kl}}_2(n\gamma , 1;q)}\, e{\lf (\f{h \gamma  }{q}\ri)} \ll (m,n,q)^{1/2}(m-n, h,q)^{1/2}q^{3/2+\varepsilon} \ene for any integers $m,n,h,q\ge1$ (see, e.g., \cite[Section 4]{KSWX}), the contribution from this case is bounded at most by a quantity
\begin{align*}
\ll &
 {X^{1+\varepsilon}\sqrt{P}}{(\T\R)^2}\sup_{\ell \ll \T\R/{\R^\ast}^{1-\varepsilon}}\,\sum_{s_1,s_2\in \ST,t \in \TT}|\beta_{s_1} \overline{\beta_{s_2}} \gamma_{t}  |^2\sum_{\substack{ r \ge 1\\ ( r ,P)=1}} (tr)^{3/2+\varepsilon}(\ell, tr)^{1/2}\\
 \ll  &{X^{1+\varepsilon}\sqrt{P}}{(\T\R)^{9/2+\varepsilon}}\mathcal{S}^2.
\end{align*}
Combining with \eqref{c42t25253155} and \eqref{nxji32w2212} shows that
\begin{align} \label{x344kvo945453} \mathcal{S}^{\text{Non-de.,1}}_{\Pi}(X,M,P) &\ll  \f{\sqrt{M}}{\mathcal{S}(\T\R)^2\R^\ast }   \Big ( X^{1+\varepsilon}\sqrt{P} \mathcal{S}\R^\ast (\T\R)^4+{X^{1+\varepsilon}\sqrt{P}}{(\T\R)^{9/2+\varepsilon}}\mathcal{S}^2\Big)^{1/2}\notag\\
&\ll  X^{1/2+\varepsilon}  \lf ( \f{\sqrt{M}P^{1/4}}{\sqrt{\mathcal{S} \R^\ast }}+\f{ \sqrt{M}(\T\R P)^{1/4}}{\R^\ast}\ri).
\end{align}
\subsubsection{Bounding $ \mathcal{S}^{\text{Non-de.,2}}_{\Pi}$}In this part, we shall turn to handling the term $ \mathcal{S}^{\text{Non-de.,2}}_{\Pi}$, the exact expression of which has been given as in \eqref{904k30fkff}. By interchanging the order of sums and applying the Cauchy-Schwarz inequality, one finds
\bee \label{x42t24c 30050} \mathcal{S}^{\text{Non-de.,2}}_{\Pi} (X,M,P)\ll \f{\sqrt{M}}{(\T\R)^{2-\varepsilon} \mathcal{S}}   \,\mathcal{B}^{1/{2}}(X,P,M),\ene
where the multiple sum $\mathcal{B}(X,P,M)$ is given by 
\begin{align}\label{9393dkdkdk} \mathcal{B}(X,P,M)=\sum_{n\ge 1} \phi {\lf ( \f{n}{\mathscr{N}} \ri)}
& \sum_{s_1,s_2\in \ST,t_1,t_2 \in \TT}\beta_{s_1}\overline{\beta_{s_2}} \gamma_{t_1}  \overline{\gamma_{t_2}}\sum_{\substack{ r_1,r_2 \ge 1\\ ( r_1r_2 ,P)=1}} \chi( r_1)  \overline{\chi(r_2 )}\notag\\
& \hskip 0.4cm   \sum_{\nu_1,\nu_2 \ge 1} \mathfrak{C}(s_1,\nu_1, n;t_1r_1) \,  \overline{\mathfrak{C}(s_2,\nu_2, n;t_2r_2) }
\,\widecheck{\mathscr{V}}\lf ( \f{\nu_1 \R^\ast }{\T\R} \ri)\notag\\
&\hskip 3.1cm 
 {\mathscr{W}}_\pm\lf ( \f{n  }{\mathscr{N}},\f{r_1}{\R}\ri) \, \overline{\,\widecheck{\mathscr{V}}\lf ( \f{\nu_2 \R^\ast }{\T\R} \ri) {\mathscr{W}}_\pm\lf ( \f{n  }{\mathscr{N}},\f{r_2}{\R}\ri)}.
\end{align}
 Recall that the co-primality relation $(t_1t_2,r_1r_2)=1$ holds. We find that there is only one scenario, namely $t_1=t_2$ and $r_1=r_2$, need to be taken into account (as far as the dominated contribution is concerned), upon recalling the explicit expression of $\mathfrak{C}$ presented in \eqref{0932jfddne2j2e2}, and employing the philosophy that $\sum_{\gamma \bmod q}{\rm{Kl}}_2(a\gamma,1;q)$ {{`vanishes'}} essentially for any integers $a,q\ge 1$. Now, if one puts $t_1=t_2=t$ and $r_1=r_2=r$, say, by applying Poisson summation with the modulus $tr$, the $n$-sum is essentially simplified into a form of the following 
 \begin{align*}&\ll \f{\mathscr{N}}{(tr)^{1-\varepsilon}} \sum_{\jmath \bmod {tr}} \mathfrak{C}(s_1,\nu_1, \jmath;tr) \, \overline{\mathfrak{C}(s_2,\nu_2, \jmath;tr) } ,
 \end{align*} upon noticing that the non-zero frequences do not exist anymore. Upon resorting to the expression \eqref{0932jfddne2j2e2}, the inner-sum over $\jmath$ turns out to be 
 \begin{align*}&\sum_{\xi, \rho\bmod {tr}}  {\rm{Kl}}_2(\overline{P}M \xi,1 ;tr)\, \overline{ {\rm{Kl}}_2(\overline{P}M \rho ,1;tr)} \, e{\lf (   \f{\nu_1 \xi-\nu_2\rho }{tr}\ri)}\sum_{\jmath \bmod {tr}}e{\lf (   \f{\jmath (\overline{s_1}\xi - \overline{s_2}  \rho  ) }{tr}\ri)}\\
 &=tr  \sum_{\xi \bmod {tr}}  {\rm{Kl}}_2(\overline{P}M\xi,1 ;tr)\,\overline{  {\rm{Kl}}_2(\overline{P}M s_2 \overline{s_1} \xi,1;tr)}\, e{\lf (   \f{(\nu_1-\nu_2 s_2 \overline{s_1}  )\xi}{tr}\ri)}.
 \end{align*}
 Next, we will proceed the argument by distinguishing two different situations.
 
 ---\emph{When $s_1=s_2$ and $\nu_1=\nu_2$.} In this case, recalling \eqref{9393dkdkdk} one finds the multiple sum $\mathcal{B}(X,P,M)$ is dominated by
 \bee \label{42cr2143by} \begin{split}&\ll (TR)^{1+\varepsilon}\mathscr{N} \,\sum_{s\in \ST,t \in \TT}|\beta_{s} \gamma_{t} |^2\sum_{\substack{ r \in \Z_{\neq  0}\\ ( r ,P)=1}} \sum_{\nu \sim \T\R/\R^\ast}  \sum_{\xi \bmod {tr}}|{\rm{Kl}}_2(\overline{P}M \xi ,1;tr)|^2
\\
 &\ll  (\T\R)^{4+\varepsilon}\mathscr{N} \mathcal{S} /\R^\ast.\end{split}
  \ene

 ---\emph{When $s_1\neq s_2$ and $\nu_1\neq \nu_2$.} Next, let us take care of the case where $s_1\neq s_2,\nu_1\neq \nu_2$; this case will provide another significant contribution to $\mathcal{B}(X,P,M)$. By appealing to the estimate \eqref{x431cr34pov2}, one infers that
  \begin{align*}\mathcal{B}(X,P,M) &\ll (TR)^{1+\varepsilon}\mathscr{N} \,\sum_{s_1,s_2\in \ST,t \in \TT}|\beta_{s_1} \overline{\beta_{s_2}} \gamma_{t}  |^2\, \sum_{\substack{ r \sim \R}} \, \sum_{\nu_1,\nu_2 \sim \T\R/\R^\ast}\,   \sum_{\xi \bmod {tr}}  \\
  &\hskip 4cm  {\rm{Kl}}_2(\overline{P}M\xi ,1;tr)\, 
  \overline{  {\rm{Kl}}_2(\overline{P}M s_2 \overline{s_1} \xi,1;tr)}\, e{\lf (   \f{(\nu_1-\nu_2 s_2 \overline{s_1}  )\xi}{tr}\ri)}\\
  &\ll  (\T\R)^{1+\varepsilon}\mathscr{N} \,\sum_{s_1,s_2\in \ST,t \in \TT}|\beta_{s_1} \overline{\beta_{s_2}} \gamma_{t}  |^2 \, \sum_{\substack{ r \sim \R}} \, \sum_{\nu_1,\nu_2 \sim \T\R/\R^\ast} (tr)^{3/2+\varepsilon} (s_1\nu _1-s_2\nu_2,tr)^{1/2}\\
  & \ll  (\T\R)^{9/2+\varepsilon}\mathscr{N} \mathcal{S}^2 /{\R^\ast}^2.
  \end{align*}Observe that, here, $s_1\nu _1\neq s_2\nu_2$; otherwise, $s_1|\nu_2$ and $s_2|\nu_1$ which cannot occur provided that $\T\R/\R^\ast<\mathcal{S}$. This is already  guaranteed by the final choice of $\R^\ast$ in \eqref{939r3c32} below. Together with \eqref{42cr2143by}, we are thus allowed to find
 \[
 \mathcal{B}(X,P,M)\ll (\T\R)^{4+\varepsilon}\mathscr{N} \mathcal{S} /\R^\ast+ (\T\R)^{9/2+\varepsilon}\mathscr{N} \mathcal{S}^2 /{\R^\ast}^2.
 \]
 This yields
 \begin{align} \label{x324114r}
 \mathcal{S}^{\text{Non-de.,2}}_{\Pi}(X,M,P) &\ll  \f{\sqrt{M}}{\mathcal{S}(\T\R)^2}   \lf (\f{(\T \R)^4 \sqrt{\mathcal{S}P}}{\sqrt{X\R^\ast}}+\f{\mathcal{S} (\T \R)^{17/4}}{\R^\ast} \sqrt{\f{P}{X}}\ri)\notag\\
 &\ll X^{1/2+\varepsilon} \mathcal{Q}^{\delta}  \lf ( \f{P\sqrt{M}}{\sqrt{\mathcal{S} \R^\ast}} +\f{P^{9/8}M^{3/4}}{\R^\ast}\ri),
 \end{align}upon combining with \eqref{x42t24c 30050}, \eqref{42cr2143by} and \eqref{234235451}. It is clear that this estimate is larger than that for $ \mathcal{S}^{\text{Non-de.,1}}_{\Pi}(X,M,P) $ as shown in \eqref{x344kvo945453}. In summary, we conclude that
 \begin{align} \label{20203ooo1}
 \mathcal{S}^{\text{Non-de.}}_{\Pi}(X,M,P)\ll X^{1/2+\varepsilon}   \mathcal{Q}^{\delta} \lf ( \f{P\sqrt{M}}{\sqrt{\mathcal{S} \R^\ast}} +\f{P^{9/8}M^{3/4}}{\R^\ast}\ri).\end{align}

  \subsection{Treatment of $ \mathcal{S}^{\text{Dua.}}_{\Pi}$.} In this subsection, we will atempt to seek a sharp estimate for $ \mathcal{S}^{\text{Dua.}}_{\Pi}$. Recall \eqref{202020}. It is straightforward to verify that \begin{align} \label{3536245564262452}  \Big| \mathcal{S}^{\text{Dua.}}_{\Pi}(X,M,P)\Big|^2\ll  &\f{1}{\T^2 M^{1-\varepsilon} } \sum_{\substack{
r_1,r_2\ge 1}} \widehat{U}\lf (\f{r_1}{M/\R}\ri) \widehat{U}\lf (\f{r_2}{M/\R}\ri)\,\sum_{\substack{ t_1,t_2\in \TT }}\gamma_{t_1} \overline{  \gamma_{t_2 }}\, 
 \notag\\
&\hskip 2.3cm \sum_{n\ge 1}\,{\rm{Kl}}_{2,\chi}(r_1,\overline{t_1}n;M)\, \overline{{\rm{Kl}}_{2,\chi}(r_2,\overline{t_2}n;M)} \,V\lf (\f{n}{X}\ri).
\end{align} By appealing to Poisson summation with the modulus $M$, the last line above is
\[ \ll \f{\widehat{V}(0)X}{M} \sum_{\delta \bmod M}{\rm{Kl}}_{2,\chi}(r_1,\overline{t_1}\delta ;M)\, \overline{{\rm{Kl}}_{2,\chi}(r_2,\overline{t_2}\delta;M)} ,\]
where $\widehat{V}$ is the Fourier transform of $V$ given by  $\widehat{V}(y)=\int_{\mathbb{R}^+}V({x})\e{{xyX}/{M}}\ud x$, and it is noticeable that the zero-frequence only survives, since the resulting length of the dual sum after Poisson is roughly $\ll_\varepsilon M^{1+\varepsilon}/X<1$. As of now,
let us take a quick glance at the exact form of the $\delta$-sum. Via opening the Kloosterman sums and executing the $\delta$-sum, the $\delta$-sum above asymptotically equals  \[\begin{split} &\chi(\overline{t_2})\,\overline{\chi}(\overline{t_1})\,M\,  \prod_{i\ge 0}\, \prod_{p^{\alpha_i}_i  \parallel  M}\,\smu_{\beta \bmod  { p^{\alpha_i}_i }}e\lf( \f{(\overline{t_1}r_1-\overline{t_2}r_2)\beta}{p^{\alpha_i}_i }\ri) \ll M^{1+\varepsilon}\, \left | (t_2r_1-t_1r_2,M) \right |. \end{split}\]  Therefore, we infer that the multiple sum on the right-hand side of \eqref{3536245564262452} is controlled by
\[ \ll X ^{1+\varepsilon}\sum_{0<r_1,r_2 \ll M/\R}\, \sum_{\substack{ t_1,t_2\in \TT\\t_2r_1\ge t_1r_2 }} (t_2r_1-t_1r_2,M).
\]
One verifies that the non-generic terms $t_1=t_2,r_1=r_2$ such that $t_2r_1=t_1r_2$ provide a quantity by $\ll_\varepsilon {X^{1+\varepsilon}M^2\T}{\R^{-1}}$; whilst, the generic terms (namely, the terms such that $t_2r_1\neq t_1r_2$), however, contribute an amount $\ll {X^{1+\varepsilon}(\T M \R^{-1})^2}$. This reveals that 
\begin{align} \label{2020345345331}
 \mathcal{S}^{\text{Dua.}}_{\Pi}(X,M,P)&\ll  X^\varepsilon   \sqrt{MX} \lf (\fone{\sqrt{\T\R}}+\f{1}{{\R}}\ri) \notag\\
&\ll \sqrt{X}\lf (\fone{p^{1/4-\varepsilon}} +\f{\T^{1+\varepsilon}}{\sqrt{MP}}\ri) , \end{align}upon recalling \eqref{234235451}.
  \subsection{Treatment of $ \mathcal{S}^{\text{Deg.}}_{\Pi}$} In this part, we are left with handling the term $ \mathcal{S}^{\text{Deg.}}_{\Pi}$ in \eqref{202010456247868}. Write $r\rightarrow rP^{\iota}$ with $\iota\ge 1$. One can recast \eqref{202010456247868} as 
\[\begin{split}
&\f{\sqrt{M}}{ (\T \R)^{1-\varepsilon}} \sum_{t\in \TT}  \gamma_t \,\sum_{\iota=1}^{\left \lfloor  \log \R/\log P \right \rfloor}\sum_{\substack{r\ge 1 }}\chi(rP^\iota) \, 
\sum_{n\ge 1} \f{{A_{\Pi}(n,1,1)}}{\sqrt{n}} \e{\f{n\overline{trP^\iota}}{M}} \U{\f{r}{
\R/P^\iota}} \V{\f{n}{X}}.\end{split}\] Here, the focal point is to analyze the situation where $\iota=1$; an entirely analogous analysis will, however, produce a much relatively smaller magnitude for the case of $\iota\ge 2$. We shall proceed by an argument which bears some resemblances to that for $ \mathcal{S}^{\text{Dua.}}_{\Pi}$ in the preceding paragraphs. In the light of this, we briefly outline the proof. Indeed, whenever $\iota=1$, one invokes the Vorono\u{\i} 
formula, Theorem \ref{9945fyhjt3464lljyrdfv432}, to the $n$-sum above. This quickly yields (essentially)  \[ \label{40404040001}\begin{split} & \f{\chi(P)}{\sqrt{M} (T R)^{1-\varepsilon}} \sum_{t\in \TT}  \gamma_t \sum_{\substack{ r \ge1}} \chi(r)\sum_{d_1|M}\,\sum_{d_2|M/d_1}\,\sum_{m\ge 1} \f{\overline{A_\Pi(m,d_2,d_1)}}{\sqrt{md^2_2d^{3}_1}}  
 \notag \\
 &\hskip 5.6cm 
{\Kl(-tr, n;d_1,d_2,{M})}\, {\mathscr{W}}_\pm\lf ( \f{m  d^2_2d^{3}_1}{\widehat{\mathscr{N}}},\f{rP}{\R}\ri)\end{split} 
 \] with the parameter $\widehat{\mathscr{N}}={M^4P}/{X^{1-\varepsilon}}$. Notice that $d_1=1$, in view of that $\widehat{\mathscr{N}}<M^3X^{-\varepsilon}$. We will concentrate on the salient scenario where
$d_2=1$, from which the dominated contribution of $ \mathcal{S}^{\text{Deg.}}_{\Pi}$ thus can be captured. Now, interchanging the order of sums, followed by the Cauchy-Schwarz inequality implies that the multiple sum above is bounded by the square-root of  
\begin {multline}  X^\varepsilon \sum_{\substack{
r_1,r_2\ge 1}}  \chi(r_1)\overline{\chi(r_1)}\sum_{\substack{
\\ t_1,t_2\in \TT }}
\gamma_{t_1} \overline{ \gamma_{t_2 }}\sum_{n \ge1 }\phi{\lf( \f{n}{\widehat{\mathscr{N}}}\ri)}  \Kl(-t_1r_1, n;{M})\\
 \overline{\Kl(-t_2r_2, n;{M})}\, {\mathscr{W}}_\pm\lf ( \f{n  }{\widehat{\mathscr{N}}}\,,\f{r_1P}{\R}\ri)\,\overline{\, {\mathscr{W}}_\pm\lf ( \f{n  }{\widehat{\mathscr{N}}}\, ,\f{r_2P}{\R}\ri)},\end{multline} 
where the weight function $\phi$ is as before. Observe that $M<\widehat{\mathscr{N}} X^{-\varepsilon}$. Equipping with  the Possion shows that the zero-frequence only exists; this implies that $t_1r_1\equiv t_2r_2 \bmod M$. If $r_1\equiv r_2 \bmod M$, it follows that the contribution to $\mathcal{S}^{\text{Deg.}}_{\Pi}(X,P,M)$ is \begin{align*}  &\ll \f{1}{\sqrt{M} (\T \R)^{1-\varepsilon}}   \lf [\f{\widehat{\mathscr{N}}}{M^{1-\varepsilon}} \sum_{\substack{
r\ll \R/P}}  \sum_{\substack{
\\ t\in \TT }}
|\gamma_{t}|^2\sum_{\gamma\bmod M} |\Kl(-tr, \gamma;{M})|^2\ri]^{1/2}\\
&\ll \sqrt{X} \mathcal{Q}^{\delta} 
P^{-3/4+\varepsilon},\end{align*} upon recalling that $(t_1r_1t_2r_2, M)=1$ and $\R,\T\ll M$ by \eqref{023e342}. On the other hand, when $r_1\neq r_2 \bmod M$, the contribution is majorized by
\begin{align}\label{jdjkd933j3ndd}  &\ll \f{1}{\sqrt{M} (\T \R)^{1-\varepsilon}}   \lf [\f{\widehat{\mathscr{N}}}{M^{1-\varepsilon}} \sum_{\substack{
r_1, r_2\ll \R/P}}  \sum_{\substack{
\\ t_1,t_2\in \TT \\t_1 \equiv t_2 r_2\overline{r_1} \bmod M}}
|\gamma_{t_1} \overline{ \gamma_{t_2 }}|\sum_{\gamma\bmod M} |\Kl(-t_1r_1, \gamma;{M})|^2\ri]^{1/2}\notag \\
&\ll \f{\sqrt{X\R} \mathcal{Q}^{\delta} }
{P^{5/4-\varepsilon}}.\end{align}
 \subsection{The endgame} At the end of this paper, we are ready to finish the proof of Theorem \ref{03}. By combining with \eqref{jdid33d33324} and the estimates in \eqref{20203ooo1},  \eqref{2020345345331} and \eqref{jdjkd933j3ndd}, one infers that
 \bee 
\label{2020432155161134441}\begin{split}\f{ \mathcal{S}_{\Pi}(X,M,P) }{X^{1/2+\varepsilon}}\ll  \,\f{ \mathcal{Q}^{\delta} P\sqrt{M}}{\sqrt{\mathcal{S} \R^\ast}} +\f{P^{9/8}M^{3/4} \mathcal{Q}^{\delta} }{\R^\ast}&+\f{\T^{1+\varepsilon}}{\sqrt{MP}}\\
 & +\fone{p^{1/4-\varepsilon}} +\f{\sqrt{\R} \mathcal{Q}^{\delta} }{P^{5/4}} \mathbf{1}_{P<\R},\end{split}\ene provided that  $\R<\mathcal{S}+\T<\R^\ast<\T \R$, and $\R+\T<M$ with $\R \T=X/M$.
A brute-force computation shows that, for the sake of securing the widest ranges of the parameters $P,M$ to ensure subconvexity, it is necessary that
\bee \label{939r3c32}\f{MP^2}{\mathcal{S}}+\mathcal{S}+\T+P^{9/8}M^{3/4}<\R^\ast < M\sqrt{P},\ene
with $\sqrt{M}P^{1/4+10\varepsilon}+\sqrt{P}<\T<\sqrt{MP} $, $\sqrt{M}<\R<\T+P^{5/2}$ and $\sqrt{M}+P^{3/2}<\mathcal{S}<\sqrt{P}M$. This implies that the choices of the parameters $\R,\R^\ast, \mathcal{S}$ and $\T$ should be
 \begin{align} \label{51515163457426583}
 \R=\f{X}{M^{3/2}P^{\mu}},\quad  \R^\ast =M^{1-\varepsilon}\sqrt{P}, \quad \mathcal{S}=M^{1-100\varepsilon}\sqrt{P},\quad \T=\sqrt{M}P^{\mu}
 \ \end{align} for any $0<\mu<1/2$, whereby the right-hand side of \eqref{2020432155161134441} turns out to be 
 \begin{align}
 \label{2r42435v6b5763}\ll X^{\varepsilon}\lf ( \mathcal{Q}^{\delta} 
 \sqrt{\f{P}{M}} +\f{ \mathcal{Q}^{\delta} P^{5/8} }{M^{1/4}}+\fone{P^{1/4}} +\fone{P^{1/2-\mu}}+\f{ \mathcal{Q}^{\delta} M^{1/4}}{P^{1+\mu/2}}\ri)
.  
\end{align}   This estimate above is trivial unless $M^{1/(4+2\mu)}<P<M^{2/5}$ for any $\mu<1/2$, as elaborated in Theorem \ref{03}.

\section{Appendix A}
\subsection{A. 1. Proof of Lemma \ref{9945fyhjt3464lljyrdfv432}} In the Appendix, we will give the proof of Lemma \ref{9945fyhjt3464lljyrdfv432}. The proof follows from the argument in \cite{HL} with a modification. We give a sketch without re-elaborating the definitions and backgrounds agian; one can find a clear blueprint of these infiormations in the $\GL(4)$ setting by comapring with that in \cite{HL} and \cite{C}. According to Corbett's notation, one assumes that $\mathtt{F}$ is a Maa\ss\ cusp form on $\GL(4)$ with trivial nebentypus, belonging to the space $L^2(\Gamma_1(N)\setminus \rm{SL}(4, \mathbb{R}))$, where the level $N\ge 1$ is a prime. 

In this appendix, we would like to take this opportunity to point out a minor error in Corbett's paper. The author is very grateful to E. Assing for many discussions upon this issue and his note \cite{AH}. Notice that there is an error occurred in \cite[p. 1392]{C} when changing the variables, which seems absent in Corbett's computations yielding the hyper-Kloosterman sums regarding the general $\GL(n)$ Vorono\u{\i} summation. Indeed, instead of the equality claimed in
the first display on the page 1392 of \cite{C}, one should has
\begin{align}\label{cc432kfff}K\ell_R(\gamma,t;\zeta, \xi)&=q^{\f{(n-2)(n-3)}{2}} \f{\prod_{i=2}^{n-2}   c^{(n-i)(n-i-1)/2}_i}{\prod_{i=3}^{n-1}d^{(i-1)(i-2)/2}_i} \cdot  \f{q^{n-2}}{\prod_{i=2}^{n-1} c^{n-i}_i}\,  {\rm{KL}}(\overline{a}\ell \gamma_0 , m; q,c,d).
 \end{align}
To see this, one finds 
\[\begin{split}  K\ell_R(\gamma,t;\zeta, \xi)&=|\xi_2 \zeta_p|^{n-2}_p |\xi_3\cdots \xi_n|^{-1}_p \psi_p(-\xi_2\xi^{-1}_3) \sum_{y_2\in \Lambda_{t_2\cdots t_{n-1}}} \cdots \sum_{y_{n-1} \in \Lambda_{t_2}}\\
&\hskip 0.5cm  \psi_p\lf (  (-1)^n   y \zeta^{-1}_p  (c_2c_3\cdots c_{n-1})^{-1} y^{-1}_2   \ri ) \, \lf ( \prod_{i=3}^{n-1} \psi_p( c_{n-i+1}   y_{i-1} y^{-1}_i)  \ri) \, \psi_p(c_2 y_{n-1}),
 \end{split}\]
via changing the variables $y_2 \rightarrow x_2\cdots x_{n-1}, y_3 \rightarrow x_2\cdots x_{n-2}, \cdots  ,y_{n-1} \rightarrow x_2$. One performs by continuing to change the variables $y_{n-1} \rightarrow  t_2 z_{n-1}, \ldots, y_2\rightarrow  t_1t_2\dots t_{n-1}z_2$, with \[
z_{n-1}\in (\Z_p/t^{-1}_2\Z_p)^\times, \quad z_{n-2}\in (\Z_p/(t_2t_3)^{-1}\Z_p)^\times, \ldots, z_{2}\in (\Z_p/(t_2t_3\cdots t_{n-1})^{-1}\Z_p)^\times,\] yielding 
\[\begin{split}  K\ell_R(\gamma,t;\zeta, \xi)&=|\xi_2 \zeta_p|^{n-2}_p |\xi_3\cdots \xi_n|^{-1}_p \psi_p(-\xi_2\xi^{-1}_3) \cdot  \prod_{i=2}^{n-2} |t_2\cdots t_{n-1}|_p  \sum_{z_{2}\in (\Z_p/(t_2t_3\cdots t_{n-1})^{-1}\Z_p)^\times} \\
& \hskip 0.7cm  \cdots\sum_{z_{n-1}\in (\Z_p/t^{-1}_2
\Z_p)^\times}  \psi_p\lf (   (-1)^n   y \zeta^{-1}_p  (c_2c_3\cdots c_{n-1})^{-1} z^{-1}_2 \cdot (t_2t_2 \cdot t_{n-1} )^{-1}  \ri ) \\
&\hskip 5.3cm  \lf ( \prod_{i=3}^{n-1} \psi_p( c_{n-i+1}  t_{n-i+2} z_{i-1} z^{-1}_i)  \ri) \, \psi_p(c_2 t_2 z_{n-1}).
 \end{split}\]
This completes the changes of variables described in the fourth display on the page 1392 of \cite{C}. Observes that the additional normalizing factor $\prod_{i=2}^{n-2} |t_2\cdots t_{n-1}|_p $ was missing in Corbett's computations. This, however, contributes a quantity
\[ \prod_{p\in R} \prod_{i=2}^{n-2} |t_2\cdots t_{n-1}|_p =\prod_{i=2}^{n-2}\lf |\f{q c_{n-i} c_{n-i-1}\cdots c_2}{d_{i+1} d_{i+2}\cdots d_{n-1} }\ri |^{i-1}=q^{\f{(n-2)(n-3)}{2}} \f{\prod_{i=2}^{n-2}   c^{(n-i)(n-i-1)/2}_i}{\prod_{i=3}^{n-1}d^{(i-1)(i-2)/2}_i}. \]
Consequently, with the expression \eqref{cc432kfff} in hand, instead of \cite[Theorem 1.1]{C}, one has
\begin{align}\label{c34245vc}&\sum_{m \neq 0} \f{A_f(m,c_1, \ldots, c_{n-1})}{|m|^{\f{n-1}{2}}}e{\lf  ( \f{am}{\ell q} \ri)} \phi_\infty(m)\prod_{p|N} \phi_p(m)\notag\\
& =q^{\f{(n-1)(n-2)}{2}}\prod_{i=2}^{n-1}|c_i|^{\f{(n-i)(n-3)}{2}}\sum_{\substack{m\neq 0\\ (m,N)=1}}  \sum_{r|N^\infty}\notag \\
&\hskip 0.8cm  \sum_{\substack{
(d_{n-1}, \ldots, d_2)=d\in \Z^{n-2} \\ d_i| \f{qc_2 \cdots c_{n-i+1} }{d_{n-1} \cdots d_{i+1}}, i=n-1, \ldots,2 }} \text{KL}(\overline{a\lambda _\ell}\cdot \ell r, m; q, c,d) \chi \lf (\f{\overline{m} qc_2\cdots c_{n-1}}{d_{n-1}\cdots d_2}\ri)^{-1} \notag \\
&\hskip 2cm  \f{A_f(d_{n-1}, \ldots, d_1, m)}{|m|^{\f{n-1}{2}} \prod_{i=2}^{n-1}  d_i^{\f{(n-3)i+2}{2}}}\, \mathcal{B}_{\pi_\infty, \phi_\infty}\lf ( \f{r m}{\lambda_\ell q^n} \prod_{i=2}^{n-1} \f{d^i_i}{c^{n-i}_i }\ri)\prod_{p|N} \mathcal{B}_{\pi_p, \Phi^{a/\ell q}_p} \lf ( \f{r m}{\lambda_\ell q^n} \prod_{i=2}^{n-1} \f{d^i_i}{c^{n-i}_i} \ri).
\end{align}
Here, $c=(c_2, \ldots, c_{n-1})\in  \Z^{n-2}$, $(a\ell,q)=(q, NM)=1$ and $\ell| (NM)^\infty $. 

Next, let us proceed by taking $M=N$, $ \phi_{\infty}(x)= |x|^{1/2}\omega{({x}/{X})} $, and setting $\phi_p=1$, \[\Phi _p(y)=\psi_p(y)W_{{{\mathtt{F}}},p}\left(\left(\begin{matrix}y & & & \\  & 1 &  &  \\ &  & 1 \\ & &  &1\end{matrix}\right)\right) \]  for $x\in \mathbb{Q}_p$ at the places $p|N$, where $W_{{{\mathtt{F}}},p}$ is a $\psi_p$-Whittaker function belonging to
$W( \Pi_{{{\mathtt{F}}},p}, \psi_p)$,  and $\omega$ is a smooth function, compactly supported on $[ 1/2 , 5/2 ]$ with bounded derivatives. In addition, we take $\chi=1$, $c_2=1$, $\ell=1$ and $q=c$ in \eqref{c34245vc}, it follows that
\begin{align}\label{4135136162626}
 &\sum_{m{\neq 0}} \f{A_{\mathtt{F}}(m,1,1)} {|m|}e\left(\frac{an}{c}\right) \omega\left(\frac{m}{X}\right)\notag \\
 &
	=c^3\sum_{\substack{m,r{\neq 0}\\ (m,N)=1\\ r\mid N^{\infty}}}\sum_{d_1\mid c}\sum_{d_2|c/d_1}
	\frac{A_{\mathtt{F}}(d_1,d_2,m)}{| m |^{3/2} d^2_2d^{5/2}_1}  \, \mathcal{KL}_2(\overline{aNL^4}r, m,c;(1,1) , (d_1,d_2) )   \notag \\ 
& \hskip 6cm \, \mathcal{B}_{\Pi_{\mathtt{F},\infty}, \phi_{\infty}}\left(\frac{rmd^2_2d^3_1}{c^4NL^4}\right)\, \prod_{p\mid N} \mathcal{B}_{\Pi_{\mathtt{F},p},\Phi^{a/c} _p}\left(\frac{rmd^2_2d^3_1}{c^4NL^4}\right),
\end{align} 
akin to \cite[Eqn. (2.9)]{HL}. Here, the function $\mathcal{B}_{\Pi_{{\mathtt{F}},\infty},\phi_{\infty}} $ is given by
\begin{align*}
	\mathcal{B}_{\Pi_{{\mathtt{F}},\infty},\phi_{\infty}}(x)= \sqrt{x}\,\mathcal{H}_{\Pi_{{\mathtt{F}},\infty}} [\omega] \left({x}{X}\right)/X
\end{align*}
with 
\[
	\mathcal{H}_{\Pi_{{\mathtt{F}},\infty}} [\omega](y) = \frac{1}{4\pi i}\sum_{\rho=0,1} \text{sgn}(y)^\rho\int_{\text{Re}( s)=\sigma}\gamma(1-s,\text{sgn}^\rho \times \Pi_{\mathtt{F},\infty},\psi_{\infty})\,\vert y\vert^{1-s}\,\widetilde{\omega}(1-s) \ud s	
\]  for any $y\in \mathbb{R}^+$, which coincides with the form $\mho_{\pm} (y;\omega)$ as shown in \S2.2 by comparing with Miller-Schmid's work \cite{MS3}. Here, $\widetilde{\omega}$ is just the normal Mellin-transform of $\omega$ given by $\widetilde{\omega}(s)=\int_{\mathbb{R}^+}\omega(x) x^{s-1}\ud x $, and $L_{\infty}(s,\text{sgn}\times \Pi_{F,\infty} )=L_{\infty}(s+1,\Pi_{F,\infty} )$. We will now need explicitly determine the expression of $p$-adic version of the Bessel transform $	\mathcal{B}_{\Pi_{{\mathtt{F}},p},\Phi^{a/c} _p}$. In fact, one finds
\begin{multline*}
	\mathcal\mathcal\mathcal{B}_{\Pi_{{\mathtt{F}},p},\Phi^{a/c} _p}(y)= \frac{\log p}{2\pi} \sum_{\substack{\xi\colon Q^{\times}_p\to C^\times \\ \xi(p)=1}}\xi(y)\int_{\sigma-i\frac{\pi}{\log p}}^{\sigma+i\frac{\pi}{\log p}}\varepsilon({1}/{2},\xi\Pi_{{\mathtt{F}},p},\psi_p)\, p^{a(\xi \Pi_{\mathtt{F},p})\cdot (s-1/{2})}\, \vert y \vert_p^{3/2-s} \\ \frac{L(s,\xi^{-1}{\Pi}_{{\mathtt{F}},p})}{L(1-s,\Pi_{{\mathtt{F}},p})}
	\,  \int_{\mathbb{Q}_p^{\times}} \xi(x)\psi_p\left (\frac{ax}{c}\right)W_{{\mathtt{F}},p}\left(\left(\begin{matrix} x & & & \\  & 1 &  &  \\ &  & 1 \\ & &  &1\end{matrix}\right)\right)\vert x\vert_p^{-1/2-s} \, d^{\times} x    d s, 
\end{multline*}
akin to \cite[Eqn. (2.12)]{HL}. Observing that $a(\xi \Pi_{{\mathtt{F}},p})=1$ in the prime level case, we are thus led to
\begin{equation}
	\mathcal\mathcal\mathcal{B}_{\Pi_{{\mathtt{F}},p},\Phi^{{a}/{c}}_{p}}(y)= \frac{\log p}{2\pi}\int_{\sigma-i\frac{\pi}{\log p}}^{\sigma+i\frac{\pi}{\log p}}p^{s-1/2}\, L(s,{\Pi}_{{\mathtt{F}},p})\, \vert y \vert_p^{3/2-s} ud s, \nonumber
\end{equation}which further shows the following exact expression that 
\begin{equation*}\begin{split}
	\mathcal\mathcal\mathcal{B}_{\Pi_{{\mathtt{F}},p},\Phi^{a/c}_{p}}(y) &= \frac{\log p}{2\pi}\sum_{k\geq 0}s_{(k,0,0,0)}({\alpha})\int_{\sigma-i\frac{\pi}{\log p}}^{\sigma+i\frac{\pi}{\log p}}p^{(1-k)s-1/2}\, \vert y \vert_p^{3/2-s} d s \\
	&= {\mathbf{1}_{v_p(y)\geq -1 } }\,{p^{-1/2}}\, \vert y \vert^{3/2}_p \, s_{(v_p(y)+1,0,0,0)}({\alpha}) \\
	&=p\,   W_{{{\mathtt{F}}},p}\left(\left(\begin{matrix} p^{v_p(y)+1} & & & \\  & 1 &  &  \\ &  & 1 \\ & &  &1\end{matrix}\right)\right).\end{split}
\end{equation*}Returning to \eqref{4135136162626}, we find that $v_{p}(rmd^2_2d^3_1/(c^4NL^4)  )=v_{p}(r/L^4)-1$, upon noticing that $(md_1d_2c,N)=1$. This yields \begin{align*} \mathcal\mathcal\mathcal{B}_{\Pi_{{\mathtt{F}},p},\Phi^{a/c} _p}\left(\frac{rmd^2_2d^3_1}{c^4NL^4}\right) =p W_{{{\mathtt{F}}},p}\left(\left(\begin{matrix} p^{v_p(r/L^4)} & & & \\  & 1 &  &  \\ &  & 1 \\ & &  &1\end{matrix}\right)\right).\end{align*}
The right-hand side of \eqref{4135136162626} thus reads\begin{align}
	& c ^3\sum_{\substack{m\neq 0\\(m,N)=1\\ L^4\mid r\mid N^{\infty}}} \sum_{d_1\mid c}\sum_{d_2|c/d_1}
	 \frac{A_{\mathtt{F}}(d_1,d_2,m)}{| m |^{3/2} d^2_2d^{5/2}_1}  \, \prod_{p\mid N}\mathcal{B}_{\Pi_{{\mathtt{F}},p},\Phi^{a/c} _p}\left(\frac{rmd^2_2d^3_1}{c^4NL^4}\right)\notag\\
	&\hskip 3cm\mathcal{KL}_2(\overline{aNL^4}r, m,c;(1,1) , (d_1,d_2) ) \,\, \mathcal{B}_{\Pi_{\mathtt{F},\infty}, \phi_{\infty}}\left(\frac{rmd^2_2d^3_1}{c^4NL^4}\right)\notag\\
	&=\frac{c^3N }{X} \sum_{\substack{m\neq 0\\(m,N)=1\\ L^4\mid r\mid N^{\infty}}} \sum_{d_1\mid c}\sum_{d_2|c/d_1}
	\frac{A_{\mathtt{F}}(d_1,d_2,m)}{| m |^{3/2} d^2_2d^{5/2}_1} \cdot \sqrt{\frac{rmd^2_2d^3_1}{c^4NL^4} } \cdot \prod_{p\mid N}   W_{{{\mathtt{F}}},p}\left(\left(\begin{matrix} p^{v_p(r/L^4)} & & & \\  & 1 &  &  \\ &  & 1 \\ & &  &1\end{matrix}\right)\right)\notag\\
	&\hskip 3cm \mathcal{KL}_2(\overline{aNL^4}r, m,c;(1,1) , (d_1,d_2) ) \,\mathcal{H}_{\Pi_{{\mathtt{F}},\infty}} [\omega]\left(\frac{rmXd^2_2d^3_1}{c^4NL^4}\right),
	\tag{4}
\end{align} which, by changing the variable $r \rightarrow rL^4$ turns out to be
\begin{align*}
	&  \frac{ c\sqrt{N}}{X} \sum_{\substack{m\neq 0\\(m,N)=1\\  r\mid N^{\infty}}} \sum_{d_1\mid c}\sum_{d_2|c/d_1}
	\frac{A_{\mathtt{F}}(d_1,d_2,m)}{| m r| d_2d_1 }  \,r^{3/2} \prod_{p\mid N}  W_{{{\mathtt{F}}},p}\left(\left(\begin{matrix}p^{v_p(r)}  & & & \\  & 1 &  &  \\ &  & 1 \\ & &  &1\end{matrix}\right)\right)\\
	&\hskip 5cm \mathcal{KL}_2(\overline{aN}, mr,c;(1,1) , (d_1,d_2) ) \,\mathcal{H}_{\Pi_{{\mathtt{F}},\infty}} [\omega]\left(\frac{rmXd^2_2d^3_1}{c^4N}\right)\\
	&=\frac{ c\sqrt{N}}{X}\sum_{\substack{m\neq 0}} \sum_{d_1\mid c}\sum_{d_2|c/d_1}
	\frac{A_{\mathtt{F}}(d_1,d_2,m)}{| m | d_2 d_1}  \,  \mathcal{KL}_2(\overline{aN}, m,c;(1,1) , (d_1,d_2) ) \\
	&\hskip 10cm \mathcal{H}_{\Pi_{{\mathtt{F}},\infty}} [\omega]\left(\frac{mXd^2_2d^3_1}{c^4N}\right).
\end{align*}
 This is essentially \begin{align*}&\frac{c\sqrt{N}}{X}\sum_{\substack{m\neq 0}} \sum_{d_1\mid c}\sum_{d_2|c/d_1}
	\frac{A_{\mathtt{F}}(d_1,d_2,m)}{| m | d_2d_1}  \,  \mathcal{KL}_2(\overline{aN}, m,c;(1,1) , (d_1,d_2) ) \notag\\
	&\hskip 8cm \mho\left ( \frac{mXd^2_2 d^3_1}{c^4N}  ;\omega\right). \end{align*}
The required assertion in Lemma \ref{9945fyhjt3464lljyrdfv432} follows immediately.

\end{document}